\documentclass[12pt]{article}
\usepackage[latin1]{inputenc}
\usepackage{amssymb,amsmath,amsfonts}
\usepackage{xcolor}
\newtheorem{theorem}{Theorem}
\newtheorem{proposition}[theorem]{Proposition}
\newtheorem{lemma}[theorem]{Lemma}

\newtheorem{remark}[theorem]{Remark}

\newcommand{\R}{\mathbb{R}}

\newcommand{\C}{\mathbb{C}}

\newcommand{\spa}{\mbox{span\,}}

\def\lp{{\langle\!\langle}}\vspace{2ex}
\def\rp{{\rangle\!\rangle}}
\def\<{{\langle}}
\def\>{{\rangle}}
\def\Sal{{\cal S}}

\def\a{\alpha}
\def\be{\begin{equation} }
\def\ee{\end{equation} }

\def\proof{\noindent\emph{Proof: }}
\def\qed{\ifhmode\unskip\nobreak\fi\ifmmode\ifinner
\else\hskip5 pt \fi\fi\hbox{\hskip5 pt \vrule width4 pt
height6 pt  depth1.5 pt \hskip 1pt }}

\begin{document}

\title{Minimal real Kaehler submanifolds}
\author{S. Chion and M. Dajczer\footnote{This research is a result of the activity developed within the framework of the Programme in Support of Excellence Groups of the Regi\'on de Murcia, Spain, by Fundaci\'on S\'eneca, Science and Technology Agency of the Regi\'on de Murcia. Marcos Dajczer was partially supported
by MICINN/FEDER project PGC2018-097046-B-I00 and Fundaci\'on S\'eneca project 19901/GERM/15, Spain. Sergio Chion was partially supported by Fundaci\'on S\'eneca project 19901/GERM/15, Spain.}}
\date{}
\maketitle

\begin{center}
\noindent
\begin{minipage}{0.85\textwidth}\parindent=15.5pt

\smallskip
\begin{center}
\large{\textsl{Dedicated to Professor Renato Tribuzy\\ on the occasion of his 75th birthday}}
\end{center}

\smallskip

{\small{
\noindent {\bf Abstract.} We show that generic rank conditions on the second fundamental 
form of an isometric immersion $f\colon M^{2n}\to\R^{2n+p}$ of a Kaehler 
manifold of complex dimension $n\geq 2$ into Euclidean space with low 
codimension $p$ imply that the submanifold has to be minimal. If 
$M^{2n}$ if simply connected, this amounts to the existence of a  
one-parameter associated family of isometric minimal immersions
unless $f$ is holomorphic.}
\smallskip

\smallskip

}

\end{minipage}
\end{center}



\section{Introduction} 

An isometric immersion $f\colon M^{2n}\to\R^{2n+p}$ is called a 
\emph{real Kaehler submanifold} if $(M^{2n},J)$ is a Kaehler 
manifold of complex dimension $n\geq 2$ immersed into 
Euclidean space with codimension $p$. We are interested in the 
case when $f$ is minimal but not holomorphic. By the latter condition 
we mean that $p$ is even and $f$ is holomorphic with respect to a 
constant complex structure in $\R^{2n+p}$.

Minimal real Kaehler submanifolds have been intensively studied since
in \cite{DG1} it was shown that they enjoy several of the basic 
properties of Euclidean minimal surfaces. For instance, if simply 
connected a minimal real Kaehler submanifold $f\colon M^{2n}\to\R^{2n+p}$ 
is either holomorphic or has a nontrivial one-parameter associated family 
of minimal isometric immersions, all of them carrying the same oriented 
Gauss map. 
Moreover, $f$ can be realized as the ``real part'' of its holomorphic 
representative $\sqrt{2}F=(f,\bar{f})\colon M^{2n}\to\C^{2n+p}$ where 
$\bar{f}\colon M^{2n}\to\R^{2n+p}$ is the conjugate immersion to $f$
in the associated family.

Real Kaehler hypersurfaces $f\colon M^{2n}\to\R^{2n+1}$, in particular
the minimal ones, have been parametrically classified by Dajczer and 
Gromoll \cite{DG1} in terms of pseudoholomorphic surfaces in spheres 
by means of the so called Gauss parametrization. A parametric 
classification of the complete minimal Real Kaehler submanifold in 
codimension two was obtained by Dajczer and Gromoll in \cite{DG2}. 
A local Weierstrass type representation for the minimal real Kaehler 
submanifolds  of any possible codimension was given by Arezzo, Pirola 
and Solci in \cite{APS}; see the Appendix of Chapter 15 in \cite{DT}.
\vspace{1ex}

We have to recall some definitions.
Let $f\colon M^{2n}\to\R^{2n+p}$ be a real Kaehler submanifold 
and let $\a\colon TM\times TM\to N_fM$ denote its second fundamental 
form taking values in the normal bundle of $f$. 
The \emph{first normal space} $N_1^\kappa(x)\subset N_fM(x)$ of $f$ at 
$x\in M^{2n}$, $\kappa\leq p$, is defined as  
$$
N_1^\kappa(x)=\spa\{\a(X,Y)\colon X,Y\in T_xM\}.
$$
If $U^s\subset N_1^\kappa(x)$ is an $s$-dimensional vector subspace
we denote
$\a_{U^s}=\pi_{U^s}\circ\a$ where $\pi_{U^s}\colon N_1^\kappa(x)\to U^s$ 
is the projection. Then let $\mathcal{N}_c(\a_{U^s})\subset T_xM$ 
be the complex vector tangent subspace given by
$$
\mathcal{N}_c(\a_{U^s})=\{Y\in T_xM\colon\a_{U^s}(X,Y)
=\a_{U^s}(X,JY)=0\;\;\mbox{for all}\;\;X\in T_xM\}
$$
and $\nu^c(\a_{U^s})=\dim\mathcal{N}_c(\a_{U^s})$.
The \emph{complex $s$-nullity} $\nu^c_s(x)$ of $f$
at $x\in M^{2n}$, $1\leq s\leq\kappa$, is defined by 
$$
\nu^c_s(x)=\max_{U^s\subset N_1(x)}\nu^c(\a_{U^s}).
$$ 
Recall that $\nu^c_\kappa(x)$ is called 
the index of complex relative nullity of $f$.

\begin{theorem}\label{min}
Let $f\colon M^{2n}\to\R^{2n+p}$, $p\leq n$, be a real Kaehler 
submanifold. Assume that at each point $x\in U$ of an open 
dense subset $U$ of $M^{2n}$ we have that either
\vspace{2ex}

\noindent $(i)$ $\nu^c_1(x)<2n-2$ and $\nu^c_{2s}(x)<2n-4s$ for
$s\geq 1$, or
\vspace{2ex}

\noindent $(ii)$ $\nu^c_{2s+1}(x)<2n-4s-2$ for $s\geq 0$.
\vspace{2ex}

\noindent Then $f$ is a minimal submanifold.
\end{theorem}

We observe that Renato Tribuzy to whom this paper is dedicated 
has given many valuable contributions to the subject of isometric 
immersions of Kaehler manifolds, for instance, see \cite{BEFT}, 
\cite{EFT}, \cite{EKT}, \cite{TF} and \cite{FT2}.

\section{Flat bilinear forms}

Let $\varphi\colon U\times V\to W$ be a bilinear form between finite 
dimensional real vector spaces. We denote by
$$
\mathcal{S}(\varphi)
=\spa\{\varphi(X,Y)\colon X\in U,\;Y\in V\},
$$
the vector subspace of $W$ generated by the \emph{image} of $\varphi$. 
The (right) \emph{kernel} of $\varphi$ is the vector subspace of 
$V$ defined by
$$
\mathcal{N}(\varphi)=\{Y\in V\colon\varphi(X,Y)=0
\;\mbox{for all}\;X\in U\}
$$
and the \emph{nullity} of $\varphi$ is
$\nu(\varphi)=\dim\mathcal{N}(\varphi)$.

A vector $X\in U$ is called a (left) \emph{regular element} 
of $\varphi$ if $\dim \varphi_X(V)=r_o$ where 
$\varphi_X\colon V\to W$ is the linear map defined by 
$\varphi_XY=\varphi(X,Y)$ and
\be\label{r}
r_o=\max\{\dim\varphi_Z(V)\colon Z\in U\}.
\ee
The set $RE(\varphi)$ of regular elements of $\varphi$ 
is easily seen to be an open dense subset of $U$; 
cf.\  Proposition $4.4$ in \cite{DT}. 
\vspace{1ex}

Let $\varphi\colon V\times V\to W$ be a bilinear form where $W$
is endowed with an inner product $\<\,,\,\>$ of any signature.
Then $\varphi$ is called \emph{flat} if 
$$
\<\varphi(X,Y),\varphi(Z,T)\>-\<\varphi(X,T),\varphi(Z,Y)\>=0
$$
for any $X,Y,Z,T\in V$.

\begin{lemma}\label{bilfla}
Let $\varphi\colon V\times V\to W$ be a flat bilinear form. If
$X\in RE(\varphi)$ then
\be\label{secondst}
\Sal(\varphi|_{V\times\ker\varphi_X})\subset\varphi_X(V)
\cap\varphi_X(V)^\perp.
\ee
\end{lemma}

\proof See equations $(8)$ and $(9)$ in \cite{Mo} or 
Proposition $4.6$ in \cite{DT}.\qed

\begin{proposition}\label{almostcomplex} 
Let $V$ be a real vector space endowed with a complex 
structure, that is, there is $J\in\text{End}(V)$ such that $J^2=-I$.
If  $X_1,JX_1,\ldots,X_{k-1},JX_{k-1},X_k\in V$ are linearly independent 
vectors then also $X_1,JX_1,\ldots,X_k,JX_k$ are linearly independent. 
In particular, we have that $V$ has even dimension.
\end{proposition}

\proof If  $JX_k=\sum_{i=1}^ka_iX_i+\sum_{j=1}^{k-1}b_jJX_j$
for $0\neq(a_1,\ldots,a_k,b_1,\ldots,b_{k-1})$ then
$$
(1+a_k^2)X_k+\sum_{j=1}^{k-1}((a_ka_j-b_j)X_j+(a_kb_j+a_j)JX_j)=0,
$$
and this is a contradiction.\vspace{2ex}\qed 

In the sequel $U^p$ denotes a $p$-dimensional vector space 
endowed with a positive definite inner product. Then 
$W^{p,p}=U^p\oplus U^p$ is endowed with the inner 
product of signature $(p,p)$ given by
$$
\lp(\xi_1,\xi_2),(\eta_1,\eta_2)\rp
=\<\xi_1,\eta_1\>_{U^p}-\<\xi_2,\eta_2\>_{U^p}.
$$

\begin{proposition}\label{even} 
Let $\rho\colon V^{2n}\times V^{2n}\to U^p$ be a symmetric  
bilinear form and $J\in\text{End}(V)$ a complex 
structure. Let $\sigma\colon V^{2n}\times V^{2n}\to W^{p,p}$ 
be the associated bilinear form given by 
$$
\sigma(X,Y)=(\rho(X,Y),\rho(X,JY))
$$  
and let $V_1,V_2\subset V^{2n}$ be vector subspaces where
$V_2$ is $J$-invariant. Then the  bilinear form 
$\sigma_0=\sigma|_{V_1\times V_2}\colon V_1\times V_2
\to\Sal(\sigma_0)\subset W^{p,p}$ satisfies:
\begin{itemize}
\item[(i)]  The vector subspace $\mathcal{N}(\sigma_0)$ of $
V_2$ is $J$-invariant.
\item [(ii)] There exists a complex structure
$T\in\text{End}\,(\Sal(\sigma_0))$.
\item [(iii)] The vector subspace $\Sal(\sigma_0)\subset W^{p,p}$    
has even dimension.
\end{itemize}
\end{proposition}

\proof We prove part $(ii)$. If $(\xi,\eta)\in\Sal(\sigma_0)$ 
let $X_i\in V_1$, $Y_i\in V_2,\; 1\leq i\leq\ell$, be such that
$$
(\xi,\eta)=\sum_{i=1}^\ell\sigma_0(X_i,Y_i)
=\sum_{i=1}^\ell(\xi_i,\eta_i).
$$
Then
$$
\sum_{i=1}^\ell\sigma_0(X_i,JY_i)
=\sum_{i=1}^\ell(\eta_i,-\xi_i)
=(\eta,-\xi)\in\Sal(\sigma_0),
$$
and hence $T\in\text{End}(\Sal(\sigma_0))$
defined by  $T(\xi,\eta)=(\eta,-\xi)$
satisfies $T^2=-I$. Now part $(iii)$ follows from
Proposition \ref{almostcomplex}.  
\vspace{2ex}\qed

Let $\a\colon V^{2n}\times V^{2n}\to U^p$ be a bilinear
form and let $J\in\text{End}(V)$ be a complex structure.
In the sequel $\beta\colon V^{2n}\times V^{2n}\to W^{p,p}$ is 
the associated bilinear form given by
\be\label{defbeta}
\beta(X,Y)=(\a(X,Y)+\a(JX,JY),\a(X,JY)-\a(JX,Y)).
\ee 
Notice that Proposition \ref{even} applies to $\beta$. It
follows from Proposition \ref{almostcomplex} that the subspace
$$
\mathcal{N}(\beta)=\{X\in V^{2n}:\a(X,JY)
=\a(JX,Y)\;\;\mbox{for all}\;\;Y\in V^{2n}\}
$$
is even dimensional.
We also have that $\beta$ verifies 
\be\label{beta0}
\beta(X,X)=(\xi,0),\;\;\beta(X,JX)=-(0,\xi),
\ee
\be\label{betasim2}
\beta(X,Y)=(\xi,\eta)\;\;\mbox{if and only if}
\;\;\beta(Y,X)=(\xi,-\eta)
\ee 
and
\be\label{betasim1}
\beta(X,Y)=(\xi,\eta)\;\;\mbox{if and only if}
\;\;\beta(X,JY)=(\eta,-\xi).
\ee
From \eqref{betasim1} we obtain that
\be\label{betasim3} 
\beta(X,Y)=\beta(S,T)\;\;\mbox{if and only if}\;\; 
\beta(X,JY)=\beta(S,JT)
\ee 
whereas from \eqref{betasim2} and \eqref{betasim1} that
\be\label{betasim4}
\beta(X,Y)=\beta(JX,JY).
\ee

\begin{proposition}\label{kerbeta}
Assume that $\beta\colon V^{2n}\times V^{2n}\to W^{p,p}$ is flat. 
If a vector subspace  $L\subset V^{2n}$ satisfies 
$\beta|_{L\times L}=0$ then $L\subset\mathcal{N}(\beta)$. 
\end{proposition}

\proof If $\beta(X,Y)=(\xi,\eta)$, we obtain from 
\eqref{betasim2} that
$$
\lp\beta(X,Y),\beta(Y,X)\rp=\|\xi\|^2_{U^p}+\|\eta\|^2_{U^p}.
$$
Thus
\be\label{condbeta}
\beta(X,Y)=0\;\;\mbox{if and only if}
\;\;\lp\beta(X,Y),\beta(Y,X)\rp=0.
\ee 
Hence, since 
$$
\lp\beta(X,Z),\beta(Z,X)\rp=\lp\beta(Z,Z),\beta(X,X)\rp=0
$$
if $Z\in L$ and $X\in V^{2p}$, then $\beta(X,Z)=0$.\qed 

\begin{proposition}\label{maincostum}
If $\beta\colon V^{2n}\times V^{2n}\to W^{p,p}$, $p\leq n$, is flat
then $\nu(\beta)=2n-r_o$ where $r_o$ is given by \eqref{r}. 
In particular, we have $\nu(\beta)\geq 2n-\dim\Sal(\beta)$.
\end{proposition}

\proof  Let $L=\ker B_X$ where $X\in RE(\beta)$ and 
$B_X\colon V^{2n}\to W^{p,p}$ is the linear map given by
$B_XY=\beta(X,Y)$. Notice that $\mathcal{N}(\beta)\subset L$. 
By \eqref{secondst} we have
$$
\lp\beta(X,Y),\beta(Y,X)\rp=0
$$
for any $X,Y\in L$. Then \eqref{condbeta} gives
$\beta|_{L\times L}=0$ and therefore $L=\mathcal{N}(\beta)$ by 
Proposition \ref{kerbeta}. Then $\nu(\beta)=\dim L=2n-r_o$.\qed 
\vspace{2ex}

The following result gives an alternative presentation and
proof of Lemma $7$ in \cite{FHZ}.

\begin{proposition}\label{diagonalization}
Assume that $\beta\colon V^{2n}\times V^{2n}\to W^{p,p}$, 
$p\leq n$, is flat. If $\nu(\beta)=2(n-p)$  there exists a basis 
$\{X_i,JX_i\}_{1\leq i\leq n}$ of $V^{2n}$ such that
\begin{itemize}
\item[(i)] $\mathcal{N}(\beta)=\spa\{X_j,JX_j,\; p+1\leq j\leq n\}$.
\item[(ii)] $\beta(Y_i,Y_j)=0\;\;\mbox{if}\;\; Y_k\in\spa\{X_k,JX_k\}
\;\mbox{where}\; i\neq j\;\mbox{and}\;k=i,j$. 
\item[(iii)] $\{\beta(X_j,X_j),\beta(X_j,JX_j),1\leq j\leq p\}$ 
is an orthonormal basis of $W^{p,p}$.
\end{itemize}
\end{proposition}

\proof Assume that the result holds for $p=n$. By part $(i)$
of Proposition \ref{even} the vector subspace $\mathcal{N}(\beta)$ is 
$J$-invariant. By Proposition \ref{almostcomplex} there is
a decomposition $V^{2n}=V_0^{2p}\oplus\mathcal{N}(\beta)$ 
where $V_0^{2p}$ is $J$-invariant. The bilinear form $\beta_0
=\beta|_{V_0\times V_0}\colon V_0^{2p}\times V_0^{2p}\to W^{p,p}$ 
is flat and $\mathcal{N}(\beta_0)=0$. In fact, if 
$Z\in\mathcal{N}(\beta_0)$ decompose 
$X\in V^{2n}$ as $X=X_1+X_2$ with 
$X_1\in V_0^{2p}$ and $X_2\in\mathcal{N}(\beta)$.
Since  $\beta(Z,X_2)=0$ we have from \eqref{betasim2} that 
$\beta(X_2,Z)=0$.  Then $Z\in \mathcal{N}(\beta)$ since
$\beta(X,Z)=\beta(X_1,Z)+\beta(X_2,Z)=0$, and thus $Z=0$.

By the initial assumption there exists a basis 
$\{X_j,JX_j\}_{1\leq j\leq p}$ of $V_0^{2p}$ such that parts 
$(ii)$ and $(iii)$ hold.
Then, by Proposition~\ref{almostcomplex} we can complete 
the basis of $V_0^{2p}$ to a basis $\{X_j,JX_j\}_{1\leq j\leq n}$ 
of $V^{2n}$ such that also part $(i)$ holds.
 
By the above, it remains to argue for the case $p=n$, 
that is, when $\nu(\beta)=0$.
\vspace{1ex}

\noindent{\it Fact} $1$.
If $p\geq 2$ there exist non-zero vectors $X,Y\in V^{2p}$ 
such that $\beta(X,Y)=0$. 
\vspace{1ex}

If $X\in\text{RE}(\beta)$ and since $\nu(\beta)=0$, then from
Proposition \ref{maincostum} the map $B_X\colon V^{2p}\to W^{p,p}$ 
is an isomorphism. Since
$\text{RE}(\beta)$ is open and dense in $V^{2p}$
there is a basis $Z_1,\ldots,Z_{2p}$ of $V^{2p}$ such that 
$Z_2\not\in\spa\{Z_1,JZ_1\}$ and  
$\{\beta(Z_k,Z_j)\}_{1\leq j\leq 2p}$ is for $k=1$ as
well as for $k=2$ a basis of $W^{p,p}$. 
Let $A=(a_{ij})$ be the $2p\times 2p$ matrix given by
$$
\beta(Z_2,Z_j)=\sum_{r=1}^{2p}a_{rj}\beta(Z_1,Z_r).
$$
Let $\lambda\in\C$ be a eigenvalue of $A$ where $(v^1,\ldots,v^{2p})\in\C^{2p}$ is the 
corresponding eigenvector.
Extending $\beta$ linearly from 
$V^{2p}\otimes\C$ to $W^{p,p}\otimes\C$, we have
$$
\sum_{j=1}^{2p}v^j\beta(Z_2,Z_j)
=\lambda\sum_{j=1}^{2p}v^j\beta(Z_1,Z_j).
$$
Hence $\beta(S,T)=0$ where $S=Z_2-\lambda Z_1$ and 
$T=\sum_{j=1}^{2p}v^jZ_j$. Then
\be\label{kernelarray}
\beta(S_1,T_1)=\beta(S_2,T_2)\;\;\mbox{and}\;\;  
\beta(S_1,T_2)+\beta(S_2,T_1)=0
\ee
where $S=S_1+iS_2$ and $T=T_1+iT_2$.
If $X=S_1-JS_2$ and $Y=T_1+JT_2$, we obtain using \eqref{betasim3},
\eqref{betasim4} and \eqref{kernelarray} that
\begin{align*}
\beta(X,Y)
&=\beta(S_1,T_1)+\beta(S_1,JT_2)-\beta(JS_2,T_1)-\beta(JS_2,JT_2)\\
&=\beta(S_1,T_1)+\beta(S_1,JT_2)+\beta(S_2,JT_1)-\beta(S_2,T_2)=0.
\end{align*}
Similarly, for $X'=S_1+JS_2$ and $Y'=T_1-JT_2$ we obtain $\beta(X',Y')=0$.  
The vectors $X$ and $X'$ are both non-zero. For instance, if $X=0$ then 
$$
JS_2+iS_2=S=Z_2-\lambda Z_1.
$$
Thus $JS_2=Z_2-\text{Re}(\lambda)Z_1$ and $S_2=-\text{Im}(\lambda)Z_1$.
Then  $Z_2\in\spa\{Z_1,JZ_1\}$, and this is a contradiction.     
Finally, if $Y=Y'=0$ then $T=0$, and this is a contradiction.  
\vspace{1ex}

\noindent{\it Fact} $2$.
There exists $Z_0\in V^{2p}$ such that $\dim\ker B_{Z_0}= 2(p-1)$. 	
\vspace{1ex}

Fact $2$ holds for $p=1$. In fact, given $0\neq X\in V^2$ we have 
from Proposition \ref{kerbeta}  that $\beta(X,X)\neq 0$. 
From \eqref{beta0} the vectors $\beta(X,X),\beta(X,JX)$ are linearly 
independent, and thus $\ker B_X=0$.

For $p\geq 2$ we argue by induction. Assume that Fact $2$ is
true for any $q\leq p-1$. By Fact $1$ there are nonzero
vectors $X,Y\in V^{2p}$ such that $\beta(X,Y)=0$. 
By part $(iii)$ of Proposition \ref{even} the dimension of 
$B_X(V)$ is even. 
If $\dim B_X(V)=2r$ then $r<p$ since $B_XY=0$. Moreover,
we have that $N_X=\ker B_X\neq V^{2p}$.  If otherwise,
we would have from \eqref{betasim2} that $X\in\mathcal{N}(\beta)$
and hence $X=0$.  Thus $\dim N_X=2p-2r$, $1\leq r\leq p-1$.

Let $U_1^s=\pi_1(B_X(V))$ where $\pi_1\colon W^{p,p}\to U^p$ 
is the projection onto the first component. We claim that $s=r$ and that
\be\label{betaximage}
B_X(V)=U_1^r\oplus U_1^r=\{\beta(Z,X)\colon Z\in V^{2p}\}.
\ee
To prove the claim, we first show that
\be\label{UsUs}
B_X(V)+\{\beta(Z,X)\colon Z\in V^{2p}\}
=U_1^s\oplus U_1^s,\;\;s\geq r.
\ee
If $\beta(X,Z)=(\xi,\eta)$ then \eqref{betasim1} gives 
$\eta\in U_1^s$, and hence $(\xi,\eta)\in U_1^s\oplus U_1^s$. 
From \eqref{betasim2} if $\beta(Z,X)=(\zeta,\eta)$ then 
$\zeta\in U_1^s$. Moreover, since $\beta(X,JZ)=-(\eta,\zeta)$
from \eqref{betasim2} and \eqref{betasim1}, then 
$\eta\in U_1^s$, and thus $(\zeta,\eta)\in U_1^s\oplus U_1^s$.  
For the other inclusion, let $(\xi_1,\xi_2)\in U_1^s\oplus U_1^s$. 
Then there are $Z_1,Z_2\in V^{2p}$ such that 
$\beta(X,Z_i)=(\xi_i,\eta_i)$, $i=1,2$. 
Then using \eqref{betasim2} and \eqref{betasim1} we obtain that
$$
(\xi_1,\xi_2)=\frac{1}{2}(\beta(X,Z_1-JZ_2)+\beta(Z_1+JZ_2,X)),
$$
and \eqref{UsUs} has been proved. If $U^p=U_1^s\oplus U_2^{p-s}$ 
is an orthogonal splitting, we show that 
\be\label{cl}
\Sal(\beta|_{N_X\times N_X})\subset U_2^{p-s}\oplus U_2^{p-s}.
\ee
The flatness of $\beta$ gives
$$
\lp\beta(X,Z),\beta(S,T)\rp=\lp\beta(X,T),\beta(S,Z)\rp=0
$$ 
for any $S,T\in N_X$ and $Z\in V^{2p}$. Moreover, from \eqref{betasim2} 
we have $\beta(S,X)=0$, and thus
$$
\lp\beta(Z,X),\beta(S,T)\rp=\lp\beta(Z,T),\beta(S,X)\rp=0,
$$
and then \eqref{UsUs} gives \eqref{cl}.

$N_X$ is $J$-invariant by Proposition \ref{even} and
$\beta|_{N_X\times N_X}\colon N_X^{2p-2r}\times N_X^{2p-2r}
\to U_2^{p-s}\oplus U_2^{p-s}$ is flat. Then Proposition \ref{maincostum} 
gives that $L=\mathcal{N}(\beta|_{N_X\times N_X})$ satisfies 
$\dim L\geq 2s-2r\geq 0$. On the other hand, since
$\beta|_{L\times L}=0$ it follows from Proposition \ref{kerbeta} that
$L\subset\mathcal{N}(\beta)=0$, and thus $s=r$. Then \eqref{betaximage} 
holds since the first equality follows from \eqref{UsUs} and 
the second equality by \eqref{betasim2}.

The assumption of induction applies to
$\beta|_{N_X\times N_X}\colon N_X\times N_X
\to U_2^{p-r}\oplus U_2^{p-r}$ since $r\geq 1$ and 
$\mathcal{N}(\beta|_{N_X\times N_X})$=0. Therefore there 
exists $Z_0\in N_X$ such that
\be\label{basis}
\dim\ker B_{Z_0}|_{N_X}=2(p-r-1).
\ee
We have that
$$
\lp\beta(S,X),\beta(Z_0,T)\rp=\lp\beta(S,T),\beta(Z_0,X)\rp=0
$$
for any $S,T\in V^{2p}$.  It follows from \eqref{betaximage} 
that $B_{Z_0}(V)\subset U_2^{p-r}\oplus U_2^{p-r}$.
Proposition \ref{kerbeta} gives $B_{Z_0}Z_0\neq 0$. Hence, 
in view of \eqref{basis} there is a basis
$X_1=Z_0,X_2=JZ_0,X_3,\ldots, X_{2(p-r)}$ of $N_X$ such that
$\beta(Z_0,X_j)=0$, $3\leq j\leq 2(p-r)$. Since 
$\mathcal{N}(\beta|_{N_X\times N_X})=0$, we
have from Proposition \ref{maincostum} that
\be\label{spanU2}
\spa\{\beta(X_i,X_j)\; 1\leq i,j\leq 2(p-r)\}
=U_2^{p-r}\oplus U_2^{p-r}.
\ee
By \eqref{beta0} we may set $\beta(Z_0,Z_0)=(\xi,0)$. 
Then \eqref{betasim2}, \eqref{betasim1} and \eqref{betasim4} give
\be\label{simspec}
\beta(Z_0,JZ_0)=(0,-\xi),\;\; \beta(JZ_0,Z_0)
=(0,\xi)\;\;\mbox{and}\;\;
\beta(JZ_0,JZ_0)=(\xi,0). 
\ee
Flatness yields
\be\label{factflat}
\lp\beta(X_s,X_t),\beta(Z_0,W)\rp
=\lp\beta(X_s,W),\beta(Z_0,X_t)\rp=0
\ee
for $3\leq s,t\leq 2(p-r)$ and any $W\in V^{2p}$.
From \eqref{spanU2}, \eqref{simspec} and \eqref{factflat} we obtain
\begin{align*}
\spa\{\beta(X_s,X_t)&\colon 3\leq s,t\leq 2(p-r)\}\\
&=(\spa\{\xi\})^\perp\cap U_2^{p-r}
\oplus(\spa\{\xi\})^\perp\cap U_2^{p-r}.
\end{align*}
It follows from \eqref{spanU2} and \eqref{factflat} that
$B_{Z_0}(V)=\spa\{\xi\}\oplus\spa\{\xi\}$, and this gives
the proof of Fact $2$. 

We conclude the proof by means of a recursive 
construction. Notice that it suffices to construct an 
orthogonal basis of $W^{p,p}$ since by \eqref{beta0}  
it can be replaced by an orthonormal one.
By Fact $2$ there is $X_1\in V^{2p}$ such that
$N_{X_1}=\ker B_{X_1}$ satisfies $\dim N_{X_1}=2p-2$ 
and by Proposition \ref{even} the vector subspace $N_{X_1}$ is $J$-invariant. 
Proposition \ref{kerbeta} gives $\beta(X_1,X_1)=(\xi_1,0)\neq 0$  and 
\eqref{betasim1} that $\beta(X_1,JX_1)=(0,-\xi_1)$.  
If $p=1$, then $X_1,JX_1$  is the desired basis. 
If $p\geq 2$ we have by flatness that
$$
\lp\beta(X_1,X_1),\beta(N_{X_1},N_{X_1})\rp=0
=\lp\beta(X_1,JX_1),\beta(N_{X_1},N_{X_1})\rp.
$$
Since $\mathcal{N}(\beta)=0$ we have from Proposition \ref{kerbeta} 
that the bilinear form 
$$
\hat\beta
=\beta|_{N_{X_1}\times N_{X_1}}\colon N_{X_1}\times N_{X_1}
\to(\spa\{\xi_1\})^\perp\oplus(\spa\{\xi_1\})^\perp
$$
satisfies $\mathcal{N}(\hat\beta)=0$.
By Fact 2 there is $X_2\in N_{X_1}$ such that 
$\dim\ker\hat{B}_{X_2}=2p-4$. As above, we have that 
$\beta(X_2,X_2)=(\xi_2,0)\neq 0$ 
and $\beta(X_2,JX_2)=(0,-\xi_2)$, where $\xi_2$ is
perpendicular to $\xi_1$. Since $N_{X_1}$ is $J$-invariant, 
then
$$
\beta(X_1,X_2)=0=\beta(X_1,JX_2).
$$
If $p=2$ then $X_1,JX_1,X_2,JX_2$ is the desired 
basis. If $p\geq 3$ we just reiterate the construction.\qed
\vspace{2ex}

In the sequel, let $\gamma\colon V^{2n}\times V^{2n}\to W^{p,p}$
be the bilinear form associated to the symmetric bilinear form
$\a\colon V^{2n}\times V^{2n}\to U^p$ given by
\be\label{defgamma}
\gamma(X,Y)=(\a(X,Y),\a(X,JY)).
\ee

\begin{proposition}\label{bgb}
Let the bilinear forms 
$\gamma,\beta\colon V^{2n}\times V^{2n}\to W^{p,p}$ be flat and
satisfy that
\be\label{productflat2}
\lp\beta(X,Y),\gamma(Z,T)\rp=\lp\beta(X,T),\gamma(Z,Y)\rp
\ee
for any $X,Y,Z,T\in V^{2n}$.  If $V_1^s=\pi_1(\Sal(\beta))$ 
and  $\a_{V_1}$ denotes taking the $V_1^s$-component of $\a$ then
$\mathcal{N}(\beta)= \mathcal{N}_c(\a_{V_1})$.
\end{proposition}

\proof We first show that
\be\label{betaimage1} 
\mathcal{S}(\beta)=V_1^s\oplus V_1^s.
\ee 
If $(\xi,\eta)\in\Sal(\beta)$ then
$(\xi,0),(0,\xi),(\eta,0)\in\Sal(\beta)$.  
In fact, if 
$$
(\xi,\eta)=\sum_k\beta(X_k,Y_k)=\sum_k(\xi_k,\eta_k),
$$
we obtain from \eqref{betasim2} and \eqref{betasim1} that
$$
\sum_k\beta(Y_k,X_k)=(\xi,-\eta),\;\;\sum_k\beta(X_k,JY_k)
=(\eta,-\xi),
$$
$$
\sum_k\beta(JY_k,X_k)
=(\eta,\xi).
$$
It follows that $\mathcal{S}(\beta)\subset V_1\oplus V_1$. 
For the other inclusion, let $(\xi,\eta)\in V_1\oplus V_1$.  
Then there are $\delta,\bar{\delta}\in U^p$ such that 
$(\xi,\delta),(\eta,\bar{\delta})\in\Sal(\beta)$, and 
by the above $(\xi,\eta)\in\Sal(\beta)$. 

From \eqref{productflat2} and \eqref{betaimage1}  we have
$\mathcal{S}(\gamma|_{V\times\mathcal{N}(\beta)})
\subset V_1^\perp\oplus V_1^\perp$.
Thus
$$
\<\a(X,Y),\xi\>=\lp\gamma(X,Y),(\xi,0)\rp=0
$$
and
$$
\<\a(X,JY),\xi\>=-\lp\gamma(X,Y),(0,\xi)\rp=0
$$
for any $X\in V^{2p}$, $Y\in\mathcal{N}(\beta)$ and $\xi\in V_1$.
Thus $\mathcal{N}(\beta)\subset\mathcal{N}_c(\a_{V_1})$.
The remaining inclusion follows from \eqref{defbeta} and 
\eqref{betaimage1}.
\qed

\section{The proof}

\begin{proposition} Let $f\colon M^{2n}\to\R^{2n+p}$ be a real 
Kaehler submanifold. 
At $x\in M^{2n}$ the bilinear forms 
$\beta,\gamma\colon T_xM\times T_xM\to W^{p,p}=N_fM(x)\oplus N_fM(x)$ 
defined by \eqref{defbeta} and \eqref{defgamma}  are 
flat and the condition \eqref{productflat2} is satisfied.
\end{proposition}

\proof The curvature tensor of a Kaehler manifold $M^{2n}$ satisfies
$$
R(X,Y)=R(JX,JY)\;\;\mbox{and}\;\;R(X,Y)JZ=JR(X,Y)Z
$$
for any $X,Y,Z\in T_xM$; cf.\ Proposition $15.1$ in \cite{DT}. 
Then straightforward computations using the Gauss equation
give the result.\qed 
\vspace{2ex}

\noindent {\em Proof of Theorem \ref{min}:} 
We claim that $V_1^s(x)=\pi_1(\Sal(\beta))$ satisfies $s=0$.
Suppose otherwise that $s>0$.  From Proposition \ref{maincostum}, 
Proposition \ref{bgb} and \eqref{betaimage1} we obtain
$$
\nu_s^c(x)\geq\nu^c(\a_{V_1^s}(x))
\geq\nu(\beta)\geq 2(n-s).
$$  
In particular, $s\neq 1$ since otherwise we have a contradiction 
with the assumptions of the theorem.

Suppose that $s\geq 2$.  If we have $\nu(\beta)>2(n-s)$ then again 
from Proposition \ref{maincostum}, 
Proposition \ref{bgb} and \eqref{betaimage1} we obtain
$$
\nu_{s-1}^c(x)\geq\nu_s^c(x)
\geq\nu^c(\a_{V_1^s}(x))\geq\nu(\beta(x))\geq 2(n-s+1)
$$
and this is contradiction with both parts of Theorem \ref{min}. 
Hence $\nu(\beta)=2(n-s)$.

\!\!Let $\{X_i,JX_i\}_{1\leq i\leq n}$ be the basis of $T_xM$ given
by Proposition \ref{diagonalization} and
$\xi_j=\pi_1(\beta(X_j,X_j))$ for $1\leq j\leq s$.  
If $i\neq j$ we obtain from \eqref{beta0} and \eqref{productflat2} that
$$
0=\lp\gamma(X,X_j),\beta(X_j,X_i)\rp=
\lp\gamma(X,X_i),\beta(X_j,X_j)\rp=\<A_{\xi_j}X_i,X\>
$$
and
$$
0=\lp\gamma(X,JX_j),\beta(X_j,X_i)\rp=\<A_{\xi_j}JX_i,X\>
$$
for any $X\in T_xM$.  Then
$$
\spa\{X_i,JX_i,\;1\leq i\leq n\;\text{and}\; i\neq j\}
\subset\ker A_{\xi_j}\cap\ker A_{\xi_j}J
$$
and thus
$$
\dim(\ker A_{\xi_j}\cap\ker A_{\xi_j}J)=2n-2,\;1\leq j\leq s.
$$
Then $\nu_{s_0}^c(x)\geq 2(n-s_0)$ for any $1\leq s_0\leq s$.  
In particular, we have $\nu_{s-1}^c(x)\geq 2(n-s+1)$ and this 
has been seen to be a contradiction.  Thus $s=0$,
that is,
\be\label{plur}
\a(JX,Y)=\a(X,JY)
\ee
for any $X,Y\in\mathfrak{X}(M)$. In particular, the 
submanifold is minimal.\qed

\begin{remark}{\em Theorem $1.2$ in \cite{DR} or Theorem $15.7$ 
in \cite{DT} give that the condition \eqref{plur} of $f$ being 
pluriharmonic is equivalent to minimality. 
}\end{remark}

\begin{remark}{\em We observe that Theorem \ref{min} does not apply
for $p=1$ since in this case the nonflat examples have $\nu_1^c=2n-2$.
 }\end{remark}


\begin{thebibliography}{99}

\bibitem{APS}
\textsc{C. Arezzo, G. Pirola  and M. Solci}, The Weierstrass representation for pluriminimal submanifolds, 
\textit{Hokkaido Math. J.} \textbf{33} (2004), 357--367. 

\bibitem{BEFT} 
\textsc{F. Burstall, J.-H. Eschenburg, M. Ferreira and R. Tribuzy}, 
K\"ahler submanifolds with parallel pluri-mean curvature,
Differential Geom. Appl. {\bf 20} (2004), 47--66.

\bibitem{DG1} 
\textsc{M. Dajczer and D. Gromoll}, Real Kaehler submanifolds and uniqueness of the Gauss map, 
\textit{J. Differential Geom.} \textbf{22} (1985), 13--28.

\bibitem{DG2} 
\textsc{M. Dajczer and D. Gromoll}, The Weierstrass representation for complete minimal real Kaehler submanifolds of codimension two, 
\textit{Invent. Math.} \textbf{119} (1995), 235--242.

\bibitem{DR} 
\textsc{M. Dajczer and L. Rodr\'{i}guez}, Rigidity of real Kaehler submanifolds,
\textit{Duke Math. J.} \textbf{53} (1986), 211--220.

\bibitem{DT} 
\textsc{M. Dajczer and R. Tojeiro}, 
\textit{``Submanifold theory beyond an introduction"},
Series: Universitext. Springer, 2019.

\bibitem{EFT} 
\textsc{J.-H. Eschenburg, M. Ferreira and R. Tribuzy}, Isotropic ppmc immersions,
\textit{Differential Geom. Appl.} \textbf{25} (2007), 351--355.

\bibitem{EKT} 
\textsc{J.-H. Eschenburg, A. Kollross and R. Tribuzy}, Codimension of immersions with parallel pluri-mean curvature, 
\textit{Differential Geom. Appl.} \textbf{27} (2009), 691--695.

\bibitem{TF}  
\textsc{M. Ferreira and R. Tribuzy}, Codimension two K\"{a}hler submanifolds of space forms,
\textit{Arch.  Math.} \textbf{79} (2002), 520--528.

\bibitem{FT2} 
\textsc{M. Ferreira and R. Tribuzy}, On the nullity of isometric immersions from K\"{a}hler manifolds, 
\textit{Rend. Semin. Mat. Univ. Politec. Torino} \textbf{65} (2007), 345--352.

\bibitem{FHZ} 
\textsc{L. Florit, W. Hui and F. Zheng}, On real Kaehler Euclidean submanifolds with non-negative Ricci curvature,
\textit{J. Eur. Math. Soc.} \textbf{7} (2005), 1--11.

\bibitem{Mo} 
\textsc{J. Moore}, Submanifolds of constant positive curvature I,
\textit{Duke Math. J.} \textbf{44} (1977), 449--489.
\end{thebibliography}
\end{document}